\documentclass[12pt]{article}
\usepackage{amssymb}
\usepackage{amsthm}
\usepackage{graphicx}
\usepackage{graphics}
\usepackage{amsthm}

\input pictex.tex
\newcommand{\esp}{\hspace{0.05cm}}

\newcommand{\f}{\bar{f}}
\newcommand{\g}{\bar{g}}
\newcommand{\clo}{\mathrm{S}^1}

\theoremstyle{definition}

\newtheorem{thm}{Theorem}[section]

\newtheorem{prop}[thm]{Proposition}

\newtheorem{lem}[thm]{Lemma}

\newtheorem{rem}[thm]{Remark}

\newtheorem{defn}[thm]{Definition}

\newcommand{\vsp}{\vspace{0.1cm}}

\begin{document}

\date{}
\author{Andr\'es Navas}

\title{A remarkable family of left-ordered groups:\\
central extensions of Hecke groups}
\maketitle

Braid groups are relevant in many branches of Mathematics.
In recent years, they have been studied as important examples of
{\em left-orderable groups} (that is, groups admitting a total order 
which is invariant by left-multiplication). Historically, the first such order on
$B_n$ (for all $n \geq 3$) was defined by Dehornoy using pure algebraic (and quite deep)
methods \cite{dehornoy-libro}. Some years later, an alternative geometric approach using 
Nielsen's theory was proposed by Thurston \cite{SW}. In this work we will, nevertheless, 
be more interested in other kinds of orders on braid groups, first
introduced by Dubrovina and Dubrovin \cite{dub}.

We will restrict the discussion to 
\esp $B_3 \!=\! \langle \sigma_1,\sigma_2 \!: 
\sigma_1 \sigma_2 \sigma_1 = \sigma_2 \sigma_1 \sigma_2 \rangle$. \esp 
(Potential generalizations to general $B_n$ will be 
discussed in \S \ref{questions}.) In \cite{dub}, it is shown 
that there is a unique left-invariant total order $\preceq_{DD}$ on $B_3$ satisfying \esp 
$\sigma_1 \sigma_2 \succ_{DD} id$ \esp and \esp $\sigma_2^{-1} \succ_{DD} id$. \esp 
This is a rather surprising fact (actually, it was conjectured as impossible in 
\cite[Conjecture 10.3.1]{pan-synt} 
which gives a new insight on the combinatorial 
structure of the Cayley graph of $B_3$ ({\em c.f.} Figure 2).

The situation described above is reminiscent to that of the Klein bottle group
$K_2 \!=\! \langle a,b \!: a^{-1} b a~=~b^{-1} \rangle$. Indeed, $K_2$ is left-orderable,
and there exits a unique left-ordering $\preceq$ satisfying $a \succ id$ and $b \succ id$.
However, $K_2$ is a less interesting example because it admits only four left-invariant
total orders (each of which is completely determined by the ``signs'' of $a$ and $b$),
whereas $B_3$ admits uncountably many \cite{braids,SW}. 

The fact that certain left-orderings are determined by finitely many inequalities
comes from the structure of their {\em positive cone}. This corresponds to the set
of elements which are {\em positive} (that is, bigger than the identity), and it
is easy to see that it is a semigroup. Actually, as is readily checked, the 
property of being left-orderable for a group $\Gamma$ is equivalent to the
existence of a (disjoint) decomposition
$$\Gamma = P^+ \sqcup P^- \sqcup \{id\},$$
where $P^+$ and $P^{-}$ are semigroups, with $P^- = \{g^{-1} \!: g \in P^+ \}$.
Now the point is that such a decomposition exists, for both $K_2$ and $B_3$, 
with $P^+$ (and $P^{-}$) finitely generated. For instance, denoting
by $\langle g_1,\ldots,g_k \rangle^+$ the semigroup generated by
$\{g_1,\ldots,g_k\}$, we have
$$K_2 = \langle a,b \rangle^+ \sqcup \langle a^{-1},b^{-1} \rangle^+ \sqcup \{id\}.$$
This decomposition can be visualized in Figure 1 below, where the
elements in $\langle a,b \rangle^+$ (that is, the positive elements
of the induced ordering) are blackened.

\vspace{0.3cm}


\beginpicture

\setcoordinatesystem units <0.67cm,0.67cm>

\putrule from -2 7 to 8 7
\putrule from -2 -1 to 8 -1
\putrule from 7 -2 to 7 8
\putrule from -1 -2 to -1 8
\putrule from 7 -2 to 7 8
\putrule from 1 -2 to 1 8
\putrule from 3 -2 to 3 2.7
\putrule from 3 3.25 to 3 8
\putrule from 2 -2 to 2 8
\putrule from 5 -2 to 5 8
\putrule from -2 1 to 8 1
\putrule from -2 3 to 2.7 3
\putrule from 3.3 3 to 8 3
\putrule from -2 4 to 8 4
\putrule from -2 5 to 8 5
\putrule from 0 -2 to 0 8
\putrule from 6 -2 to 6 8
\putrule from -2  6 to 8 6
\putrule from 4 -2 to 4 8
\putrule from -2 2 to 8 2
\putrule from 0 -2 to 0 8
\putrule from -2 0 to 8 0

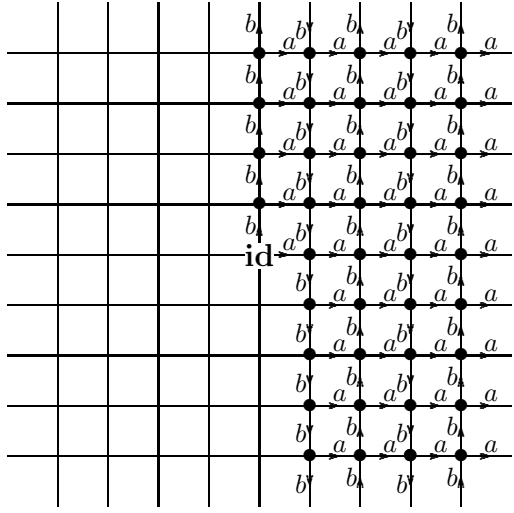
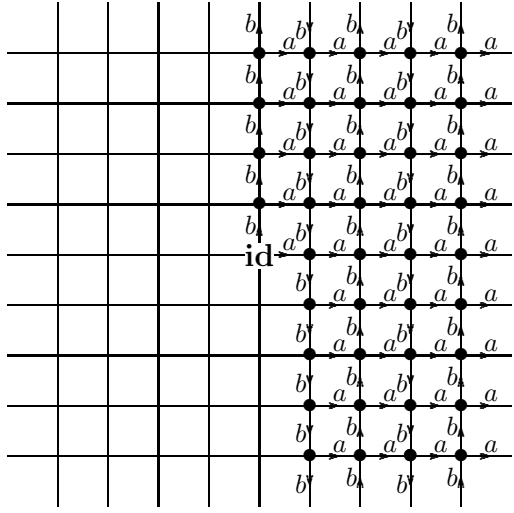
\put{Figure 1: The finitely generated positive cone $P^+ \! = \langle a,b \rangle^+$ 
on $K_2 = \langle a,b \!: a^{-1}ba = b^{-1} \rangle$} at 3 -2.8
\put{$\bf{id}$} at 3 3
\put{} at -9.4 0


\plot 3.55 7 3.4 7.05 /
\plot 3.55 7 3.4 6.95 /

\plot 3.55 6 3.4 6.05 /
\plot 3.55 6 3.4 5.95 /

\plot 3.55 5 3.4 5.05 /
\plot 3.55 5 3.4 4.95 /

\plot 3.55 4 3.4 4.05 /
\plot 3.55 4 3.4 3.95 /

\plot 3.55 3 3.4 3.05 /
\plot 3.55 3 3.4 2.95 /


\plot 4.55 7 4.4 7.05 /
\plot 4.55 7 4.4 6.95 /

\plot 4.55 6 4.4 6.05 /
\plot 4.55 6 4.4 5.95 /

\plot 4.55 5 4.4 5.05 /
\plot 4.55 5 4.4 4.95 /

\plot 4.55 4 4.4 4.05 /
\plot 4.55 4 4.4 3.95 /

\plot 4.55 3 4.4 3.05 /
\plot 4.55 3 4.4 2.95 /

\plot 4.55 2 4.4 2.05 /
\plot 4.55 2 4.4 1.95 /

\plot 4.55 1 4.4 1.05 /
\plot 4.55 1 4.4 0.95 /

\plot 4.55 0 4.4 0.05 /
\plot 4.55 0 4.4 -0.05 /

\plot 4.55 -1 4.4 -0.95 /
\plot 4.55 -1 4.4 -1.05 /


\plot 5.55 7 5.4 7.05 /
\plot 5.55 7 5.4 6.95 /

\plot 5.55 6 5.4 6.05 /
\plot 5.55 6 5.4 5.95 /

\plot 5.55 5 5.4 5.05 /
\plot 5.55 5 5.4 4.95 /

\plot 5.55 4 5.4 4.05 /
\plot 5.55 4 5.4 3.95 /

\plot 5.55 3 5.4 3.05 /
\plot 5.55 3 5.4 2.95 /

\plot 5.55 2 5.4 2.05 /
\plot 5.55 2 5.4 1.95 /

\plot 5.55 1 5.4 1.05 /
\plot 5.55 1 5.4 0.95 /

\plot 5.55 0 5.4 0.05 /
\plot 5.55 0 5.4 -0.05 /

\plot 5.55 -1 5.4 -0.95 /
\plot 5.55 -1 5.4 -1.05 /


\plot 6.55 7 6.4 7.05 /
\plot 6.55 7 6.4 6.95 /

\plot 6.55 6 6.4 6.05 /
\plot 6.55 6 6.4 5.95 /

\plot 6.55 5 6.4 5.05 /
\plot 6.55 5 6.4 4.95 /

\plot 6.55 4 6.4 4.05 /
\plot 6.55 4 6.4 3.95 /

\plot 6.55 3 6.4 3.05 /
\plot 6.55 3 6.4 2.95 /

\plot 6.55 2 6.4 2.05 /
\plot 6.55 2 6.4 1.95 /

\plot 6.55 1 6.4 1.05 /
\plot 6.55 1 6.4 0.95 /

\plot 6.55 0 6.4 0.05 /
\plot 6.55 0 6.4 -0.05 /

\plot 6.55 -1 6.4 -0.95 /
\plot 6.55 -1 6.4 -1.05 /


\plot 7.55 7 7.4 7.05 /
\plot 7.55 7 7.4 6.95 /

\plot 7.55 6 7.4 6.05 /
\plot 7.55 6 7.4 5.95 /

\plot 7.55 5 7.4 5.05 /
\plot 7.55 5 7.4 4.95 /

\plot 7.55 4 7.4 4.05 /
\plot 7.55 4 7.4 3.95 /

\plot 7.55 3 7.4 3.05 /
\plot 7.55 3 7.4 2.95 /

\plot 7.55 2 7.4 2.05 /
\plot 7.55 2 7.4 1.95 /

\plot 7.55 1 7.4 1.05 /
\plot 7.55 1 7.4 0.95 /

\plot 7.55 0 7.4 0.05 /
\plot 7.55 0 7.4 -0.05 /

\plot 7.55 -1 7.4 -0.95 /
\plot 7.55 -1 7.4 -1.05 /


\plot 4 7.45 4.05 7.6 /
\plot 4 7.45 3.95 7.6 /

\plot 6 7.45 6.05 7.6 /
\plot 6 7.45 5.95 7.6 /

\plot 4 6.45 4.05 6.6 /
\plot 4 6.45 3.95 6.6 /

\plot 6 6.45 6.05 6.6 /
\plot 6 6.45 5.95 6.6 /

\plot 4 5.45 4.05 5.6 /
\plot 4 5.45 3.95 5.6 /

\plot 6 5.45 6.05 5.6 /
\plot 6 5.45 5.95 5.6 /

\plot 4 4.45 4.05 4.6 /
\plot 4 4.45 3.95 4.6 /

\plot 6 4.45 6.05 4.6 /
\plot 6 4.45 5.95 4.6 /

\plot 4 3.45 4.05 3.6 /
\plot 4 3.45 3.95 3.6 /

\plot 6 3.45 6.05 3.6 /
\plot 6 3.45 5.95 3.6 /

\plot 4 2.45 4.05 2.6 /
\plot 4 2.45 3.95 2.6 /

\plot 6 2.45 6.05 2.6 /
\plot 6 2.45 5.95 2.6 /

\plot 4 1.45 4.05 1.6 /
\plot 4 1.45 3.95 1.6 /

\plot 6 1.45 6.05 1.6 /
\plot 6 1.45 5.95 1.6 /

\plot 4 0.45 4.05 0.6 /
\plot 4 0.45 3.95 0.6 /

\plot 6 0.45 6.05 0.6 /
\plot 6 0.45 5.95 0.6 /

\plot 4 -0.55 4.05 -0.4 /
\plot 4 -0.55 3.95 -0.4 /

\plot 6 -0.55 6.05 -0.4 /
\plot 6 -0.55 5.95 -0.4 /

\plot 4 -1.55 4.05 -1.4 /
\plot 4 -1.55 3.95 -1.4 /

\plot 6 -1.55 6.05 -1.4 /
\plot 6 -1.55 5.95 -1.4 /


\plot 3 3.55
3.05 3.4 /

\plot 3 3.55
2.95 3.4 /

\plot 3 4.55
3.05 4.4 /

\plot 3 4.55
2.95 4.4 /


\plot 3 5.55
3.05 5.4 /

\plot 3 5.55
2.95 5.4 /


\plot 3 6.55
3.05 6.4 /

\plot 3 6.55
2.95 6.4 /


\plot 5 3.55
5.05 3.4 /

\plot 5 3.55
4.95 3.4 /


\plot 5 4.55
5.05 4.4 /

\plot 5 4.55
4.95 4.4 /


\plot 5 5.55
5.05 5.4 /

\plot 5 5.55
4.95 5.4 /


\plot 5 6.55
5.05 6.4 /

\plot 5 6.55
4.95 6.4 /


\plot 5 2.55
5.05 2.4 /

\plot 5 2.55
4.95 2.4 /


\plot 5 1.55
5.05 1.4 /

\plot 5 1.55
4.95 1.4 /


\plot 5 0.55
5.05 0.4 /

\plot 5 0.55
4.95 0.4 /


\plot 5 -0.45
5.05 -0.6 /

\plot 5 -0.45
4.95 -0.6 /


\plot 7 3.55
7.05 3.4 /

\plot 7 3.55
6.95 3.4 /


\plot 7 4.55
7.05 4.4 /

\plot 7 4.55
6.95 4.4 /


\plot 7 5.55
7.05 5.4 /

\plot 7 5.55
6.95 5.4 /


\plot 7 6.55
7.05 6.4 /

\plot 7 6.55
6.95 6.4 /


\plot 7 2.55
7.05 2.4 /

\plot 7 2.55
6.95 2.4 /


\plot 7 1.55
7.05 1.4 /

\plot 7 1.55
6.95 1.4 /


\plot 7 0.55
7.05 0.4 /

\plot 7 0.55
6.95 0.4 /


\plot 7 -0.45
7.05 -0.6 /

\plot 7 -0.45
6.95 -0.6 /


\plot 5 7.55
5.05 7.4 /

\plot 5 7.55
4.95 7.4 /


\plot 7 7.55
7.05 7.4 /

\plot 7 7.55
6.95 7.4 /


\plot 3 7.55
3.05 7.4 /

\plot 3 7.55
2.95 7.4 /


\plot 5 -1.45
5.05 -1.6 /

\plot 5 -1.45
4.95 -1.6 /


\plot 7 -1.45
7.05 -1.6 /

\plot 7 -1.45
6.95 -1.6 /


\put{$\bullet$} at 3 6
\put{$\bullet$} at 3 4
\put{$\bullet$} at 3 5

\put{$\bullet$} at 4 0
\put{$\bullet$} at 4 1
\put{$\bullet$} at 4 2
\put{$\bullet$} at 4 3
\put{$\bullet$} at 4 4
\put{$\bullet$} at 4 5
\put{$\bullet$} at 4 6

\put{$\bullet$} at 5 0
\put{$\bullet$} at 5 1
\put{$\bullet$} at 5 2
\put{$\bullet$} at 5 3
\put{$\bullet$} at 5 4
\put{$\bullet$} at 5 5
\put{$\bullet$} at 5 6

\put{$\bullet$} at 6 0
\put{$\bullet$} at 6 1
\put{$\bullet$} at 6 2
\put{$\bullet$} at 6 3
\put{$\bullet$} at 6 4
\put{$\bullet$} at 6 5
\put{$\bullet$} at 6 6

\put{$\bullet$} at 7 0
\put{$\bullet$} at 7 1
\put{$\bullet$} at 7 2
\put{$\bullet$} at 7 3
\put{$\bullet$} at 7 4
\put{$\bullet$} at 7 5
\put{$\bullet$} at 7 6

\put{$\bullet$} at 3 7
\put{$\bullet$} at 4 7
\put{$\bullet$} at 5 7
\put{$\bullet$} at 6 7
\put{$\bullet$} at 7 -1
\put{$\bullet$} at 7 7
\put{$\bullet$} at 4 -1
\put{$\bullet$} at 5 -1
\put{$\bullet$} at 6 -1
\put{$\bullet$} at 6 -1


\begin{footnotesize}

\put{$a$} at 3.42 7.2
\put{$a$} at 3.42 6.2
\put{$a$} at 3.42 5.2
\put{$a$} at 3.42 4.2
\put{$a$} at 3.42 3.2

\put{$a$} at 4.42 7.2
\put{$a$} at 4.42 6.2
\put{$a$} at 4.42 5.2
\put{$a$} at 4.42 4.2
\put{$a$} at 4.42 3.2
\put{$a$} at 4.42 2.2
\put{$a$} at 4.42 1.2
\put{$a$} at 4.42 0.2
\put{$a$} at 4.42 -0.8

\put{$a$} at 5.42 7.2
\put{$a$} at 5.42 6.2
\put{$a$} at 5.42 5.2
\put{$a$} at 5.42 4.2
\put{$a$} at 5.42 3.2
\put{$a$} at 5.42 2.2
\put{$a$} at 5.42 1.2
\put{$a$} at 5.42 0.2
\put{$a$} at 5.42 -0.8

\put{$a$} at 6.42 7.2
\put{$a$} at 6.42 6.2
\put{$a$} at 6.42 5.2
\put{$a$} at 6.42 4.2
\put{$a$} at 6.42 3.2
\put{$a$} at 6.42 2.2
\put{$a$} at 6.42 1.2
\put{$a$} at 6.42 0.2
\put{$a$} at 6.42 -0.8

\put{$a$} at 7.42 7.2
\put{$a$} at 7.42 6.2
\put{$a$} at 7.42 5.2
\put{$a$} at 7.42 4.2
\put{$a$} at 7.42 3.2
\put{$a$} at 7.42 2.2
\put{$a$} at 7.42 1.2
\put{$a$} at 7.42 0.2
\put{$a$} at 7.42 -0.8


\put{$b$} at 3.65 7.46
\put{$b$} at 3.65 6.46
\put{$b$} at 3.65 5.46
\put{$b$} at 3.65 4.46
\put{$b$} at 3.65 3.46
\put{$b$} at 3.65 2.46
\put{$b$} at 3.65 1.46
\put{$b$} at 3.65 0.46
\put{$b$} at 3.65 -0.54
\put{$b$} at 3.65 -1.54

\put{$b$} at 5.65 7.46
\put{$b$} at 5.65 6.46
\put{$b$} at 5.65 5.46
\put{$b$} at 5.65 4.46
\put{$b$} at 5.65 3.46
\put{$b$} at 5.65 2.46
\put{$b$} at 5.65 1.46
\put{$b$} at 5.65 0.46
\put{$b$} at 5.65 -0.54
\put{$b$} at 5.65 -1.54


\put{$b$} at 2.65 7.6
\put{$b$} at 2.65 6.6
\put{$b$} at 2.65 5.6
\put{$b$} at 2.65 4.6
\put{$b$} at 2.65 3.6

\put{$b$} at 4.65 7.6
\put{$b$} at 4.65 6.6
\put{$b$} at 4.65 5.6
\put{$b$} at 4.65 4.6
\put{$b$} at 4.65 3.6
\put{$b$} at 4.65 2.6
\put{$b$} at 4.65 1.6
\put{$b$} at 4.65 0.6
\put{$b$} at 4.65 -0.4
\put{$b$} at 4.65 -1.4

\put{$b$} at 6.65 7.6
\put{$b$} at 6.65 6.6
\put{$b$} at 6.65 5.6
\put{$b$} at 6.65 4.6
\put{$b$} at 6.65 3.6
\put{$b$} at 6.65 2.6
\put{$b$} at 6.65 1.6
\put{$b$} at 6.65 0.6
\put{$b$} at 6.65 -0.4
\put{$b$} at 6.65 -1.4

\end{footnotesize}

\endpicture


\vspace{0.7cm}

A similar phenomenon occurs for $B_3$. Indeed, letting $a = \sigma_1 \sigma_2$ and
$b= \sigma_2^{-1}$, we also have the decomposition
$$B_3 = \langle a,b \rangle^+ \sqcup \langle a^{-1},b^{-1} \rangle^+ \sqcup \{ id \}.$$
The proof of this fact is given in \cite{dub}. It is very indirect and uses Dehornoy's theory.
We will propose an alternative argument which applies to a larger family of groups. As a
byproduct, we will retrieve (and generalize) the Dehornoy ordering and some of its
properties by rather elementary methods (see \S \ref{dehornoy}).

As in the case of $K_2$, the decomposition of $B_3$ above may be
easily illustrated: see Figure 2.
The Cayley graph of $B_3$ is, essentially, a product of $\mathbb{Z}^2$ by a dyadic rooted tree.
The (quasi-isometric) copy of $\mathbb{Z}^2$ corresponds to the ``upper level'' of the graph,
and the corresponding edges are slightly blackened. An arrow pointing to the right should be
added to every horizontal edge of the graph. These edges represent multiplications by
$a$, and all other (oriented) edges represent  multiplications by $b$.
Starting at $id$, every blackened element can be reached by a path
that follows the direction of the arrows. Conversely, every element which is not blackened 
may be reached by a path starting at $id$ following a direction opposite to that
of the arrows. Finally, no (nontrivial) element can be reached both ways.

\vspace{0.2cm}

The Klein bottle group may be presented in the form
$$K_2 = \langle a,b \! : bab = a \rangle.$$
Moreover, with respect to the generators $a = \sigma_1\sigma_2$
and $b = \sigma_2^{-1}$, the standard presentation of $B_3$ becomes
$$B_3 = \langle a,b \!: ba^2b = a \rangle.$$
This makes natural the study of the groups
$$\Gamma_n = \langle a,b\!: ba^n b = a \rangle.$$
These groups have been already considered in \cite{deh,paris} as examples lying
on the border of the theory of Gaussian and Garside groups. We will show that,
though they do not fit in these important categories, they share a remarkable 
combinatorial property with $B_3$.

\vspace{0.35cm}

\noindent{\bf Main Theorem.} {\em For each $n \geq 1$, 
the group $\Gamma_n$ admits the decomposition}
$$\Gamma_n = \langle a,b \rangle^+ \sqcup 
\langle a^{-1},b^{-1} \rangle^+ \sqcup \{ id \}.$$

\vspace{0.5cm}


\beginpicture

\setcoordinatesystem units <0.82cm,0.82cm>

\putrule from 0 0 to 11 0
\putrule from 0 1 to 11 1
\putrule from 2 0 to 2 1
\putrule from 5 0 to 5 1
\putrule from 8 0 to 8 1
\putrule from 11 0 to 11 1

\putrule from 0 1.025 to 11 1.025
\putrule from 0 0.975 to 11 0.975

\putrule from 2 -2.025 to 13 -2.025
\putrule from 2 -1.975 to 13 -1.975

\putrule from 4 -5.025 to 9.7 -5.025
\putrule from 10.3 -5.025 to 15 -5.025

\putrule from 4 -4.975 to 9.7 -4.975
\putrule from 10.3 -4.975 to 15 -4.975

\putrule from 6 -8.025 to 17 -8.025
\putrule from 6 -7.975 to 17 -7.975

\putrule from 8 -11.025 to 19 -11.025
\putrule from 8 -10.975 to 19 -10.975

\putrule from 0.6 -1.1 to 11.6 -1.1

\plot 2 0 2.6 -1.1 /
\plot 5 0 5.6 -1.1 /
\plot 8 0 8.6 -1.1 /
\plot 11 0 11.6 -1.1 /

\plot 1 0  0 1 /
\plot 4 0  3 1 /
\plot 7 0  6 1 /
\plot 10 0 9 1 /

\plot 0 0  1.6 -1.1 /
\plot 3 0  4.6 -1.1 /
\plot 6 0  7.6 -1.1 /
\plot 9 0 10.6 -1.1 /


\putrule from 2 -3 to 13 -3
\putrule from 2 -2 to 13 -2
\putrule from 4 -3 to 4 -2
\putrule from 7 -3 to 7 -2
\putrule from 10 -3 to 10 -2
\putrule from 13 -3 to 13 -2

\putrule from 2.6 -4.1 to 13.6 -4.1

\plot 4 -3 4.6 -4.1 /
\plot 7 -3 7.6 -4.1 /
\plot 10 -3 10.6 -4.1 /
\plot 13 -3 13.6 -4.1 /

\plot 3 -3  2 -2 /
\plot 6 -3  5 -2 /
\plot 9 -3  8 -2 /
\plot 12 -3 11 -2 /

\plot 2 -3  3.6 -4.1 /
\plot 5 -3  6.6 -4.1 /
\plot 8 -3 9.6 -4.1 /
\plot 11 -3 12.6 -4.1 /


\putrule from 4 -6 to 15 -6
\putrule from 4 -5 to 9.7 -5
\putrule from 10.3 -5 to 15 -5
\putrule from 6 -6 to 6 -5
\putrule from 9 -6 to 9 -5
\putrule from 12 -6 to 12 -5
\putrule from 15 -6 to 15 -5

\putrule from 4.6 -7.1 to 15.6 -7.1

\plot 6 -6 6.6 -7.1 /
\plot 9 -6 9.6 -7.1 /
\plot 12 -6 12.6 -7.1 /
\plot 15 -6 15.6 -7.1 /

\plot 5 -6  4 -5 /
\plot 8 -6  7 -5 /
\plot 11 -6  10.2 -5.2 /
\plot 14 -6 13 -5 /

\plot 4 -6  5.6 -7.1 /
\plot 7 -6  8.6 -7.1 /
\plot 10 -6 11.6 -7.1 /
\plot 13 -6 14.6 -7.1 /


\putrule from 6 -9 to 17 -9
\putrule from 6 -8 to 17 -8
\putrule from 8 -9 to 8 -8
\putrule from 11 -9 to 11 -8
\putrule from 14 -9 to 14 -8
\putrule from 17 -9 to 17 -8

\putrule from 6.6 -10.1 to 17.6 -10.1

\plot 8  -9 8.6  -10.1 /
\plot 11 -9 11.6 -10.1 /
\plot 14 -9 14.6 -10.1 /
\plot 17 -9 17.6 -10.1 /

\plot 7  -9 6 -8 /
\plot 10 -9 9 -8 /
\plot 13 -9 12 -8 /
\plot 16 -9 15 -8 /

\plot 6 -9  7.6  -10.1 /
\plot 9 -9  10.6 -10.1 /
\plot 12 -9 13.6 -10.1 /
\plot 15 -9 16.6 -10.1 /


\putrule from 8 -12 to 19 -12
\putrule from 8 -11 to 19 -11
\putrule from 10 -12 to 10 -11
\putrule from 13 -12 to 13 -11
\putrule from 16 -12 to 16 -11
\putrule from 19 -12 to 19 -11

\putrule from 8.6 -13.1 to 19.6 -13.1

\plot 10 -12 10.6 -13.1 /
\plot 13 -12 13.6 -13.1 /
\plot 16 -12 16.6 -13.1 /
\plot 19 -12 19.6 -13.1 /

\plot 9 -12  8 -11 /
\plot 12 -12  11 -11 /
\plot 15 -12 14 -11 /
\plot 18 -12 17 -11 /

\plot 8 -12  9.6 -13.1 /
\plot 11 -12 12.6 -13.1 /
\plot 14 -12 15.6 -13.1 /
\plot 17 -12 18.6 -13.1 /


\plot 1 1 9 -11 / \plot 1.025 1 9.025 -11 / \plot 0.975 1 8.975 -11 /
\plot 4 1 12 -11 / \plot 4.025 1 12.025 -11 / \plot 3.975 1 11.975 -11 /
\plot 7 1 15 -11 / \plot 7.025 1 15.025 -11 / \plot 6.975 1 14.975 -11 /
\plot 10 1 18 -11 / \plot 10.025 1 18.025 -11 / \plot 9.975 1 17.975 -11 /

\plot 3 1 4 -2 / \plot 3.025 1 4.025 -2 / \plot 2.975 1 3.975 -2 /
\plot 6 1 7 -2 / \plot 6.025 1 7.025 -2 / \plot 5.975 1 6.975 -2 /
\plot 9 1 10 -2 / \plot 9.025 1 10.025 -2 / \plot 8.975 1 9.975 -2 /

\plot 5 -2 6 -5 / \plot 5.025 -2 6.025 -5 / \plot 4.975 -2 5.975 -5 /
\plot 8 -2 9 -5 / \plot 8.025 -2 9.025 -5 / \plot 7.975 -2 8.975 -5 /
\plot 11 -2 12 -5 / \plot 11.025 -2 12.025 -5 / \plot 10.975 -2 11.975 -5 /

\plot 7 -5 8 -8 / \plot 7.025 -5 8.025 -8 / \plot 6.975 -5 7.975 -8 /
\plot 10.1 -5.2 11 -8 / \plot 10.125 -5.25 11.02 -8 / \plot 10.075 -5.2 10.975 -8 /
\plot 13 -5 14 -8 / \plot 13.025 -5 14.025 -8 / \plot 12.975 -5 13.975 -8 /

\plot 9 -8 10 -11 / \plot 9.025 -8 10.025 -11 / \plot 8.975 -8 9.975 -11 /
\plot 12 -8 13 -11 / \plot 12.025 -8 13.025 -11 / \plot 11.975 -8 12.975 -11 /
\plot 15 -8 16 -11 / \plot 15.025 -8 16.025 -11 / \plot 14.975 -8 15.975 -11 /


\plot 2.7 -1.5 2.95 -1.7 / \plot 2.7 -1.5 2.7 -1.8 /
\plot 5.7 -1.5 5.95 -1.7 / \plot 5.7 -1.5 5.7 -1.8 /
\plot 8.7 -1.5 8.95 -1.7 / \plot 8.7 -1.5 8.7 -1.8 /
\plot 11.7 -1.5 11.95 -1.7 / \plot 11.7 -1.5 11.7 -1.8 /

\plot 4.7 -4.5 4.95 -4.7 / \plot 4.7 -4.5 4.7 -4.8 /
\plot 7.7 -4.5 7.95 -4.7 / \plot 7.7 -4.5 7.7 -4.8 /
\plot 10.7 -4.5 10.95 -4.7 / \plot 10.7 -4.5 10.7 -4.8 /
\plot 13.7 -4.5 13.95 -4.7 / \plot 13.7 -4.5 13.7 -4.8 /

\plot 6.7 -7.5 6.95 -7.7 / \plot 6.7 -7.5 6.7 -7.8 /
\plot 9.7 -7.5 9.95 -7.7 / \plot 9.7 -7.5 9.7 -7.8 /
\plot 12.7 -7.5 12.95 -7.7 / \plot 12.7 -7.5 12.7 -7.8 /
\plot 15.7 -7.5 15.95 -7.7 / \plot 15.7 -7.5 15.7 -7.8 /

\plot 8.7 -10.5 8.95 -10.7 / \plot 8.7 -10.5 8.7 -10.8 /
\plot 11.7 -10.5 11.95 -10.7 / \plot 11.7 -10.5 11.7 -10.8 /
\plot 14.7 -10.5 14.95 -10.7 / \plot 14.7 -10.5 14.7 -10.8 /
\plot 17.7 -10.5 17.95 -10.7 / \plot 17.7 -10.5 17.7 -10.8 /

\plot 2 0.4 1.9 0.6 / \plot 2 0.4 2.1 0.6 /
\plot 5 0.4 4.9 0.6 / \plot 5 0.4 5.1 0.6 /
\plot 8 0.4 7.9 0.6 / \plot 8 0.4 8.1 0.6 /
\plot 11 0.4 10.9 0.6 / \plot 11 0.4 11.1 0.6 /

\plot 4 -2.6 3.9 -2.4 / \plot 4 -2.6 4.1 -2.4 /
\plot 7 -2.6 6.9 -2.4 / \plot 7 -2.6 7.1 -2.4 /
\plot 10 -2.6 9.9 -2.4 / \plot 10 -2.6 10.1 -2.4 /
\plot 13 -2.6 12.9 -2.4 / \plot 13 -2.6 13.1 -2.4 /

\plot 6 -5.6 5.9 -5.4 / \plot 6 -5.6 6.1 -5.4 /
\plot 9 -5.6 8.9 -5.4 / \plot 9 -5.6 9.1 -5.4 /
\plot 12 -5.6 11.9 -5.4 / \plot 12 -5.6 12.1 -5.4 /
\plot 15 -5.6 14.9 -5.4 / \plot 15 -5.6 15.1 -5.4 /

\plot 8 -8.6 7.9 -8.4 / \plot 8 -8.6 8.1 -8.4 /
\plot 11 -8.6 10.9 -8.4 / \plot 11 -8.6 11.1 -8.4 /
\plot 14 -8.6 13.9 -8.4 / \plot 14 -8.6 14.1 -8.4 /
\plot 17 -8.6 16.9 -8.4 / \plot 17 -8.6 17.1 -8.4 /

\plot 10 -11.6 9.9 -11.4 / \plot 10 -11.6 10.1 -11.4 /
\plot 13 -11.6 12.9 -11.4 / \plot 13 -11.6 13.1 -11.4 /
\plot 16 -11.6 15.9 -11.4 / \plot 16 -11.6 16.1 -11.4 /
\plot 19 -11.6 18.9 -11.4 / \plot 19 -11.6 19.1 -11.4 /

\plot 3.87 -1.6 3.95 -1.4 / \plot 3.87 -1.6 3.67 -1.42 /
\plot 6.87 -1.6 6.95 -1.4 / \plot 6.87 -1.6 6.67 -1.42 /
\plot 9.87 -1.6 9.95 -1.4 / \plot 9.87 -1.6 9.67 -1.42 /

\plot 5.87 -4.6 5.95 -4.4 / \plot 5.87 -4.6 5.67 -4.42 /
\plot 8.87 -4.6 8.95 -4.4 / \plot 8.87 -4.6 8.67 -4.42 /
\plot 11.87 -4.6 11.95 -4.4 / \plot 11.87 -4.6 11.67 -4.42 /

\plot 7.87 -7.6 7.95 -7.4 / \plot 7.87 -7.6 7.67 -7.42 /
\plot 10.87 -7.6 10.95 -7.4 / \plot 10.87 -7.6 10.67 -7.42 /
\plot 13.87 -7.6 13.95 -7.4 / \plot 13.87 -7.6 13.67 -7.42 /

\plot 9.87 -10.6 9.95 -10.4 / \plot 9.87 -10.6 9.67 -10.42 /
\plot 12.87 -10.6 12.95 -10.4 / \plot 12.87 -10.6 12.67 -10.42 /
\plot 15.87 -10.6 15.95 -10.4 / \plot 15.87 -10.6 15.67 -10.42 /

\plot 0.4 0.6 0.78 0.4 / \plot 0.4 0.6 0.55 0.25 /
\plot 3.4 0.6 3.78 0.4 / \plot 3.4 0.6 3.55 0.25 /
\plot 6.4 0.6 6.78 0.4 / \plot 6.4 0.6 6.55 0.25 /
\plot 9.4 0.6 9.78 0.4 / \plot 9.4 0.6 9.55 0.25 /

\plot 2.4 -2.4 2.78 -2.6 / \plot 2.4 -2.4 2.55 -2.75 /
\plot 5.4 -2.4 5.78 -2.6 / \plot 5.4 -2.4 5.55 -2.75 /
\plot 8.4 -2.4 8.78 -2.6 / \plot 8.4 -2.4 8.55 -2.75 /
\plot 11.4 -2.4 11.78 -2.6 / \plot 11.4 -2.4 11.55 -2.75 /

\plot 4.4 -5.4 4.78 -5.6 / \plot 4.4 -5.4 4.55 -5.75 /
\plot 7.4 -5.4 7.78 -5.6 / \plot 7.4 -5.4 7.55 -5.75 /
\plot 10.4 -5.4 10.78 -5.6 / \plot 10.4 -5.4 10.55 -5.75 /
\plot 13.4 -5.4 13.78 -5.6 / \plot 13.4 -5.4 13.55 -5.75 /

\plot 6.4 -8.4 6.78 -8.6 / \plot 6.4 -8.4 6.55 -8.75 /
\plot 9.4 -8.4 9.78 -8.6 / \plot 9.4 -8.4 9.55 -8.75 /
\plot 12.4 -8.4 12.78 -8.6 / \plot 12.4 -8.4 12.55 -8.75 /
\plot 15.4 -8.4 15.78 -8.6 / \plot 15.4 -8.4 15.55 -8.75 /

\plot 8.4 -11.4 8.78 -11.6 / \plot 8.4 -11.4 8.55 -11.75 /
\plot 11.4 -11.4 11.78 -11.6 / \plot 11.4 -11.4 11.55 -11.75 /
\plot 14.4 -11.4 14.78 -11.6 / \plot 14.4 -11.4 14.55 -11.75 /
\plot 17.4 -11.4 17.78 -11.6 / \plot 17.4 -11.4 17.55 -11.75 /

\plot 2.37 -0.7 2.4 -0.45 / \plot 2.37 -0.7 2.12 -0.5 /
\plot 5.37 -0.7 5.4 -0.45 / \plot 5.37 -0.7 5.12 -0.5 /
\plot 8.37 -0.7 8.4 -0.45 / \plot 8.37 -0.7 8.12 -0.5 /
\plot 11.37 -0.7 11.4 -0.45 / \plot 11.37 -0.7 11.12 -0.5 /

\plot 4.37 -3.7 4.4 -3.45 / \plot 4.37 -3.7 4.12 -3.5 /
\plot 7.37 -3.7 7.4 -3.45 / \plot 7.37 -3.7 7.12 -3.5 /
\plot 10.37 -3.7 10.4 -3.45 / \plot 10.37 -3.7 10.12 -3.5 /
\plot 13.37 -3.7 13.4 -3.45 / \plot 13.37 -3.7 13.12 -3.5 /

\plot 6.37 -6.7 6.4 -6.45 / \plot 6.37 -6.7 6.12 -6.5 /
\plot 9.37 -6.7 9.4 -6.45 / \plot 9.37 -6.7 9.12 -6.5 /
\plot 12.37 -6.7 12.4 -6.45 / \plot 12.37 -6.7 12.12 -6.5 /
\plot 15.37 -6.7 15.4 -6.45 / \plot 15.37 -6.7 15.12 -6.5 /

\plot 8.37 -9.7 8.4 -9.45 / \plot 8.37 -9.7 8.12 -9.5 /
\plot 11.37 -9.7 11.4 -9.45 / \plot 11.37 -9.7 11.12 -9.5 /
\plot 14.37 -9.7 14.4 -9.45 / \plot 14.37 -9.7 14.12 -9.5 /
\plot 17.37 -9.7 17.4 -9.45 / \plot 17.37 -9.7 17.12 -9.5 /

\plot 10.37 -12.7 10.4 -12.45 / \plot 10.37 -12.7 10.12 -12.5 /
\plot 13.37 -12.7 13.4 -12.45 / \plot 13.37 -12.7 13.12 -12.5 /
\plot 16.37 -12.7 16.4 -12.45 / \plot 16.37 -12.7 16.12 -12.5 /
\plot 19.37 -12.7 19.4 -12.45 / \plot 19.37 -12.7 19.12 -12.5 /

\plot 1.05 -0.7 1.4 -0.8 / \plot 1.05 -0.7 1.15 -0.95 /
\plot 4.05 -0.7 4.4 -0.8 / \plot 4.05 -0.7 4.15 -0.95 /
\plot 7.05 -0.7 7.4 -0.8 / \plot 7.05 -0.7 7.15 -0.95 /
\plot 10.05 -0.7 10.4 -0.8 / \plot 10.05 -0.7 10.15 -0.95 /

\plot 3.05 -3.7 3.4 -3.8 / \plot 3.05 -3.7 3.15 -3.95 /
\plot 6.05 -3.7 6.4 -3.8 / \plot 6.05 -3.7 6.15 -3.95 /
\plot 9.05 -3.7 9.4 -3.8 / \plot 9.05 -3.7 9.15 -3.95 /
\plot 12.05 -3.7 12.4 -3.8 / \plot 12.05 -3.7 12.15 -3.95 /

\plot 5.05 -6.7 5.4 -6.8 / \plot 5.05 -6.7 5.15 -6.95 /
\plot 8.05 -6.7 8.4 -6.8 / \plot 8.05 -6.7 8.15 -6.95 /
\plot 11.05 -6.7 11.4 -6.8 / \plot 11.05 -6.7 11.15 -6.95 /
\plot 14.05 -6.7 14.4 -6.8 / \plot 14.05 -6.7 14.15 -6.95 /

\plot 7.05 -9.7 7.4 -9.8 / \plot 7.05 -9.7 7.15 -9.95 /
\plot 10.05 -9.7 10.4 -9.8 / \plot 10.05 -9.7 10.15 -9.95 /
\plot 13.05 -9.7 13.4 -9.8 / \plot 13.05 -9.7 13.15 -9.95 /
\plot 16.05 -9.7 16.4 -9.8 / \plot 16.05 -9.7 16.15 -9.95 /

\plot 9.05 -12.7 9.4 -12.8 / \plot 9.05 -12.7 9.15 -12.95 /
\plot 12.05 -12.7 12.4 -12.8 / \plot 12.05 -12.7 12.15 -12.95 /
\plot 15.05 -12.7 15.4 -12.8 / \plot 15.05 -12.7 15.15 -12.95 /
\plot 18.05 -12.7 18.4 -12.8 / \plot 18.05 -12.7 18.15 -12.95 /

\plot 1 -0.15 1.05 -0.38 / \plot 1 -0.15 0.8 -0.3 /
\plot 4 -0.15 4.05 -0.38 / \plot 4 -0.15 3.8 -0.3 /
\plot 7 -0.15 7.05 -0.38 / \plot 7 -0.15 6.8 -0.3 /
\plot 10 -0.15 10.05 -0.38 / \plot 10 -0.15 9.8 -0.3 /

\plot 3 -3.15 3.05 -3.38 / \plot 3 -3.15 2.8 -3.3 /
\plot 6 -3.15 6.05 -3.38 / \plot 6 -3.15 5.8 -3.3 /
\plot 9 -3.15 9.05 -3.38 / \plot 9 -3.15 8.8 -3.3 /
\plot 12 -3.15 12.05 -3.38 / \plot 12 -3.15 11.8 -3.3 /

\plot 5 -6.15 5.05 -6.38 / \plot 5 -6.15 4.8 -6.3 /
\plot 8 -6.15 8.05 -6.38 / \plot 8 -6.15 7.8 -6.3 /
\plot 11 -6.15 11.05 -6.38 / \plot 11 -6.15 10.8 -6.3 /
\plot 14 -6.15 14.05 -6.38 / \plot 14 -6.15 13.8 -6.3 /

\plot 7 -9.15 7.05 -9.38 / \plot 7 -9.15 6.8 -9.3 /
\plot 10 -9.15 10.05 -9.38 / \plot 10 -9.15 9.8 -9.3 /
\plot 13 -9.15 13.05 -9.38 / \plot 13 -9.15 12.8 -9.3 /
\plot 16 -9.15 16.05 -9.38 / \plot 16 -9.15 15.8 -9.3 /

\plot 9 -12.15 9.05 -12.38 / \plot 9 -12.15 8.8 -12.3 /
\plot 12 -12.15 12.05 -12.38 / \plot 12 -12.15 11.8 -12.3 /
\plot 15 -12.15 15.05 -12.38 / \plot 15 -12.15 14.8 -12.3 /
\plot 18 -12.15 18.05 -12.38 / \plot 18 -12.15 17.8 -12.3 /

\plot 2.6 -0.25 2.65 -0.08 / \plot 2.6 -0.25 2.825 -0.22 /
\plot 5.6 -0.25 5.65 -0.08 / \plot 5.6 -0.25 5.825 -0.22 /
\plot 8.6 -0.25 8.65 -0.08 / \plot 8.6 -0.25 8.825 -0.22 /

\plot 4.6 -3.25 4.65 -3.08 / \plot 4.6 -3.25 4.825 -3.22 /
\plot 7.6 -3.25 7.65 -3.08 / \plot 7.6 -3.25 7.825 -3.22 /
\plot 10.6 -3.25 10.65 -3.08 / \plot 10.6 -3.25 10.825 -3.22 /

\plot 6.6 -6.25 6.65 -6.08 / \plot 6.6 -6.25 6.825 -6.22 /
\plot 9.6 -6.25 9.65 -6.08 / \plot 9.6 -6.25 9.825 -6.22 /
\plot 12.6 -6.25 12.65 -6.08 / \plot 12.6 -6.25 12.825 -6.22 /

\plot 8.6 -9.25 8.65 -9.08 / \plot 8.6 -9.25 8.825 -9.22 /
\plot 11.6 -9.25 11.65 -9.08 / \plot 11.6 -9.25 11.825 -9.22 /
\plot 14.6 -9.25 14.65 -9.08 / \plot 14.6 -9.25 14.825 -9.22 /

\plot 10.6 -12.25 10.65 -12.08 / \plot 10.6 -12.25 10.825 -12.22 /
\plot 13.6 -12.25 13.65 -12.08 / \plot 13.6 -12.25 13.825 -12.22 /
\plot 16.6 -12.25 16.65 -12.08 / \plot 16.6 -12.25 16.825 -12.22 /


\setdots


\putrule from -0.3 -0.8 to 10.7 -0.8

\plot 3 0  1.7 -0.8 /
\plot 6 0  4.7 -0.8 /
\plot 9 0  7.7 -0.8 /

\plot 1 0 0.7 -0.8 /
\plot 4 0 3.7 -0.8 /
\plot 7 0 6.7 -0.8 /
\plot 10 0 9.7 -0.8 /


\putrule from 1.7 -3.8 to 12.7 -3.8

\plot 5 -3  3.7 -3.8 /
\plot 8 -3  6.7 -3.8 /
\plot 11 -3  9.7 -3.8 /

\plot 3 -3 2.7 -3.8 /
\plot 6 -3 5.7 -3.8 /
\plot 9 -3 8.7 -3.8 /
\plot 12 -3 11.7 -3.8 /


\putrule from 3.7 -6.8 to 14.7 -6.8

\plot 7 -6  5.7 -6.8 /
\plot 10 -6  8.7 -6.8 /
\plot 13 -6  11.7 -6.8 /

\plot 5 -6 4.7 -6.8 /
\plot 8 -6 7.7 -6.8 /
\plot 11 -6 10.7 -6.8 /
\plot 14 -6 13.7 -6.8 /


\putrule from 5.7 -9.8 to 16.7 -9.8

\plot 9 -9  7.7 -9.8 /
\plot 12 -9  10.7 -9.8 /
\plot 15 -9  13.7 -9.8 /

\plot 7 -9 6.7 -9.8 /
\plot 10 -9 9.7 -9.8 /
\plot 13 -9 12.7 -9.8 /
\plot 16 -9 15.7 -9.8 /


\putrule from 7.7 -12.8 to 18.7 -12.8
\plot 11 -12  9.7 -12.8 /
\plot 14 -12  12.7 -12.8 /
\plot 17 -12  15.7 -12.8 /

\plot 9  -12 8.7  -12.8 /
\plot 12 -12 11.7 -12.8 /
\plot 15 -12 14.7 -12.8 /
\plot 18 -12 17.7 -12.8 /

\begin{small}

\put{$\bullet$} at 12.85 -11
\put{$\bullet$} at 13.9 -11
\put{$\bullet$} at 14.85 -11
\put{$\bullet$} at 15.85 -11
\put{$\bullet$} at 16.85 -11
\put{$\bullet$} at 17.82 -11
\put{$\bullet$} at 18.85 -11

\put{$\bullet$} at 12.9 -12
\put{$\bullet$} at 13.9 -12
\put{$\bullet$} at 14.9 -12
\put{$\bullet$} at 15.9 -12
\put{$\bullet$} at 16.9 -12
\put{$\bullet$} at 17.9 -12
\put{$\bullet$} at 18.9 -12

\put{$\bullet$} at 13.45 -13.1
\put{$\bullet$} at 14.5 -13.1
\put{$\bullet$} at 15.45 -13.1
\put{$\bullet$} at 16.45 -13.1
\put{$\bullet$} at 17.45 -13.1
\put{$\bullet$} at 18.45 -13.1
\put{$\bullet$} at 19.45 -13.1

\put{$\bullet$} at 12.6 -12.8
\put{$\bullet$} at 13.6 -12.8
\put{$\bullet$} at 14.6 -12.8
\put{$\bullet$} at 15.6 -12.8
\put{$\bullet$} at 16.6 -12.8
\put{$\bullet$} at 17.6 -12.8
\put{$\bullet$} at 18.6 -12.8

\put{$\bullet$} at 10.88 -8
\put{$\bullet$} at 11.9 -8
\put{$\bullet$} at 12.85 -8
\put{$\bullet$} at 13.85 -8
\put{$\bullet$} at 14.85 -8
\put{$\bullet$} at 15.85 -8
\put{$\bullet$} at 16.9 -8

\put{$\bullet$} at 10.9 -9
\put{$\bullet$} at 11.9 -9
\put{$\bullet$} at 12.9 -9
\put{$\bullet$} at 13.9 -9
\put{$\bullet$} at 14.9 -9
\put{$\bullet$} at 15.9 -9
\put{$\bullet$} at 16.9 -9

\put{$\bullet$} at 11.45 -10.1
\put{$\bullet$} at 12.5 -10.1
\put{$\bullet$} at 13.45 -10.1
\put{$\bullet$} at 14.45 -10.1
\put{$\bullet$} at 15.4 -10.1
\put{$\bullet$} at 16.5 -10.1
\put{$\bullet$} at 17.45 -10.1

\put{$\bullet$} at 10.6 -9.8
\put{$\bullet$} at 11.6 -9.8
\put{$\bullet$} at 12.6 -9.8
\put{$\bullet$} at 13.6 -9.8
\put{$\bullet$} at 14.6 -9.8
\put{$\bullet$} at 15.6 -9.8
\put{$\bullet$} at 16.6 -9.8

\put{$\bullet$} at 10.84 -5
\put{$\bullet$} at 11.85 -5
\put{$\bullet$} at 12.85 -5
\put{$\bullet$} at 13.85 -5
\put{$\bullet$} at 14.9 -5

\put{$\bullet$} at 11.9 -6
\put{$\bullet$} at 12.9 -6
\put{$\bullet$} at 13.9 -6
\put{$\bullet$} at 14.9 -6

\put{$\bullet$} at 12.45 -7.1
\put{$\bullet$} at 13.34 -7.1
\put{$\bullet$} at 14.4 -7.1
\put{$\bullet$} at 15.45 -7.1

\put{$\bullet$} at 11.6 -6.8
\put{$\bullet$} at 12.6 -6.8
\put{$\bullet$} at 13.6 -6.8
\put{$\bullet$} at 14.6 -6.8

\put{$\bullet$} at 8.84 -2
\put{$\bullet$} at 9.85 -2
\put{$\bullet$} at 10.85 -2
\put{$\bullet$} at 11.85 -2
\put{$\bullet$} at 12.9 -2

\put{$\bullet$} at 9.9 -3
\put{$\bullet$} at 10.9 -3
\put{$\bullet$} at 11.9 -3
\put{$\bullet$} at 12.9 -3

\put{$\bullet$} at 9.6 -3.8
\put{$\bullet$} at 10.6 -3.8
\put{$\bullet$} at 11.6 -3.8
\put{$\bullet$} at 12.6 -3.8

\put{$\bullet$} at 10.45 -4.1
\put{$\bullet$} at 11.34 -4.1
\put{$\bullet$} at 12.4 -4.1
\put{$\bullet$} at 13.45 -4.1

\put{$\bullet$} at 6.87 1
\put{$\bullet$} at 7.85 1
\put{$\bullet$} at 8.85 1
\put{$\bullet$} at 9.85 1
\put{$\bullet$} at 10.9 1

\put{$\bullet$} at 7.9 0
\put{$\bullet$} at 8.9 0
\put{$\bullet$} at 9.9 0
\put{$\bullet$} at 10.9 0

\put{$\bullet$} at 8.45 -1.1
\put{$\bullet$} at 9.34 -1.1
\put{$\bullet$} at 10.4 -1.1
\put{$\bullet$} at 11.45 -1.1

\put{$\bullet$} at 7.6 -0.8
\put{$\bullet$} at 8.6 -0.8
\put{$\bullet$} at 9.6 -0.8
\put{$\bullet$} at 10.6 -0.8

\end{small}
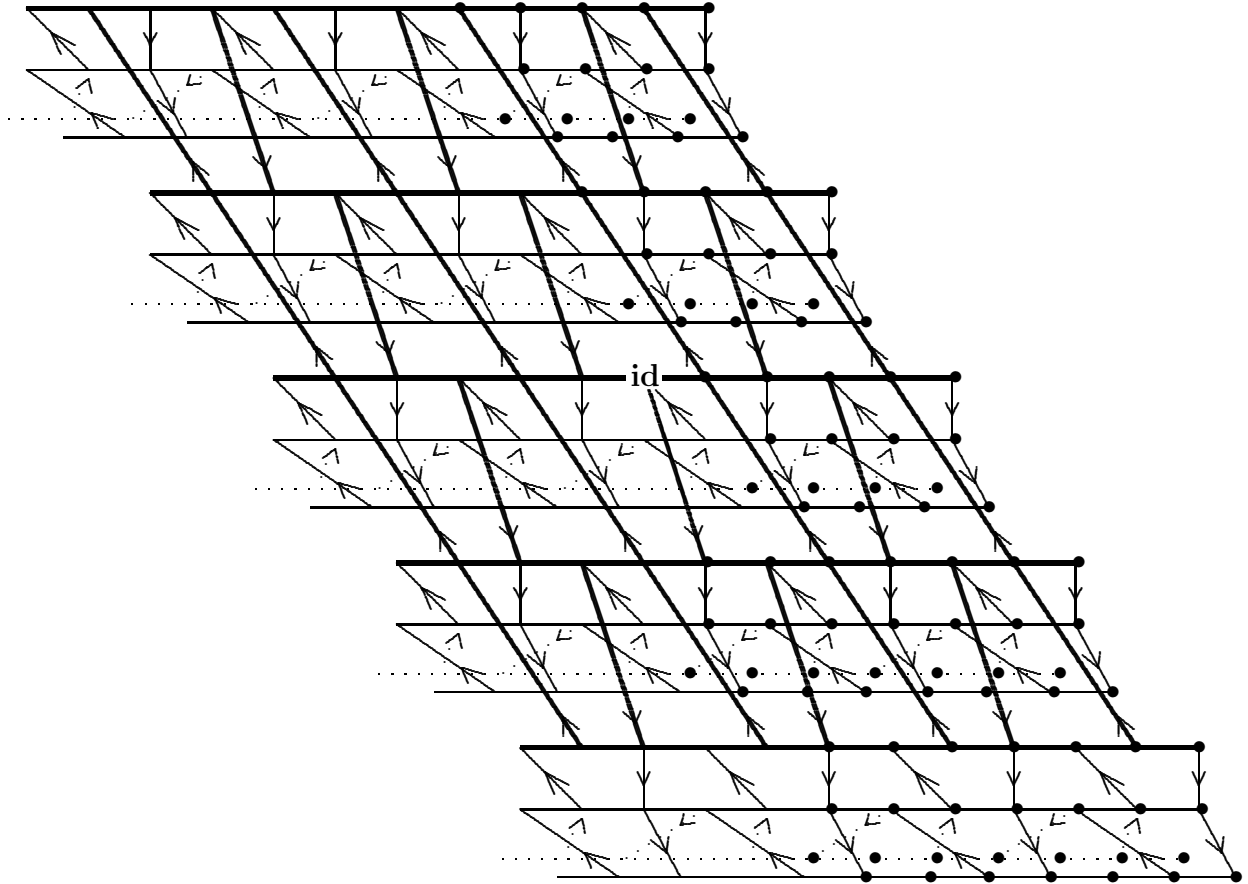
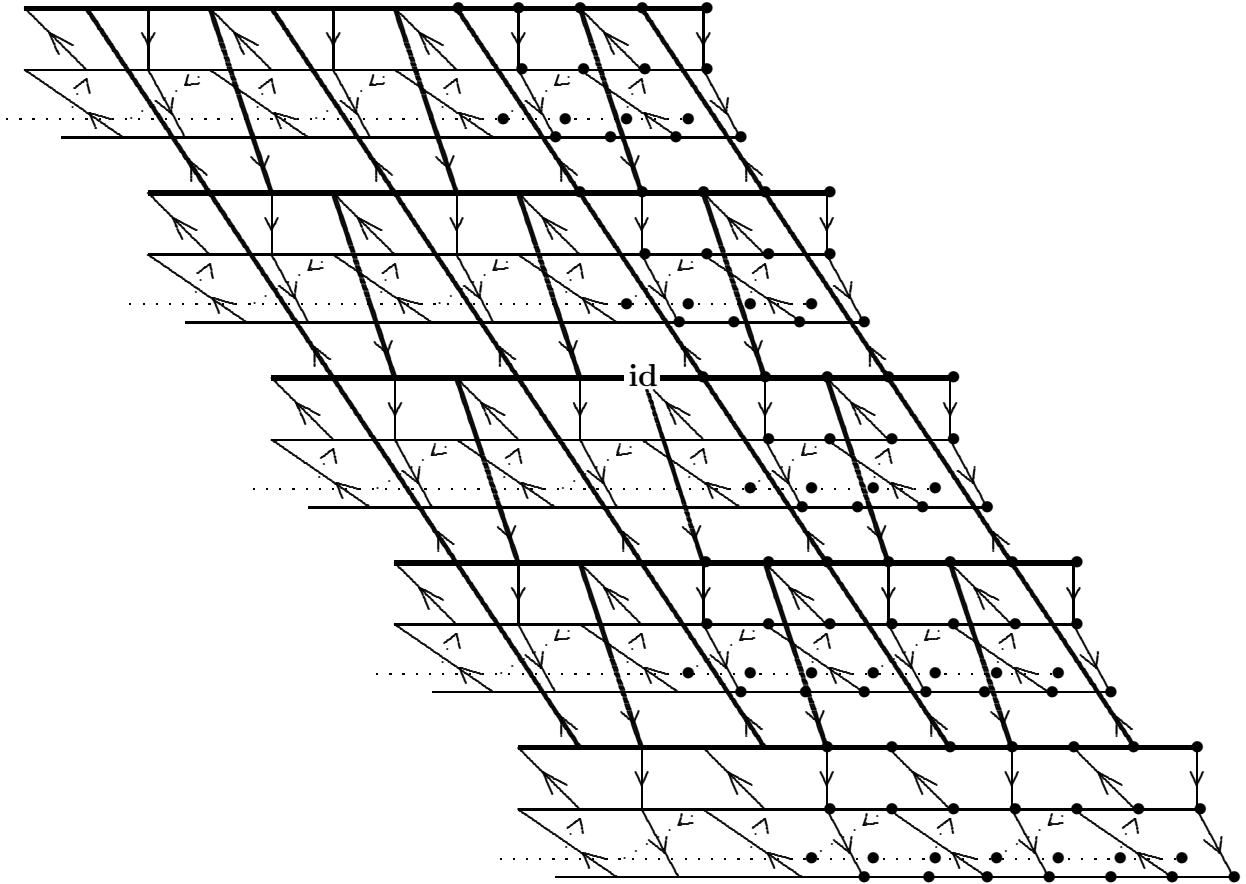
\put{Figure 2: The Cayley graph of
$B_3 = \langle a,b \!: b a^2 b = a \rangle$ and the DD-positive cone}
at 10 -14

\put{{\bf id}} at 9.7 -4.95

\endpicture


\vspace{0.65cm}

The proof of this result involves two issues. First, we need to show that every nontrivial element
$w \!\in\! \Gamma_n$ belongs to either $\langle a,b \rangle^+$ or $\langle a^{-1},b^{-1} \rangle^+$.
For this, we begin by appealing to the theory of Garside groups and write $w$ in the form
$w = uv^{-1}$ for some $u,v$ in $\langle a,b \rangle^+$ (see \S \ref{ore}). This creates
central {\em handles}, that is, expressions of the form $ab^ka^{-1}$. The main point
here is that these handles belong to one of the semigroups above. Indeed, from the
relation $b a^n b = a$ one easily deduces that $a b^k a^{-1} = (a^{-(n-1)}b^{-1})^{k}$. Due to
this, the corresponding {\em reduction procedure} that we will perform in \S \ref{en-uno-solo} for showing
that \esp $\Gamma_n \setminus \{id\} = \langle a,b \rangle^+ \cup \langle a^{-1},b^{-1} \rangle^+$, \esp
though similar, is much simpler than the Dehornoy {\em handle reduction algorithm} \cite{DDRW}. In
particular, its convergence follows from elementary combinatorial arguments.

The second issue consists in showing that $\langle a,b \rangle^+$ and $\langle a^{-1},b^{-1} \rangle^+$
are disjoint, which is equivalent to showing that $\langle a,b \rangle^+$ does not contain the
identity. This is done in \S \ref{positive} by means of a quite simple ping-pong type argument.
As a motivation, recall the well-known representation of $B_3$ in
$\mathrm{PSL}(2,\mathbb{R})$ given by
$$\sigma_1 \longrightarrow
\left[
\begin{array}
{cc}
1 & 1  \\
0 & 1  \\
\end{array}
\right],\qquad
\sigma_2 \longrightarrow
\left[
\begin{array}
{cc}
1 & 0  \\
-1 & 1  \\
\end{array}
\right].$$
This representation induces an action of $B_3$ on the circle (viewed as the projective
line), and by looking at the dynamics of this action, this yields the desired
property by a ping-pong type argument.
The extension of this proof to $\Gamma_n$ is straightforward, as one may
easily produce an action of $\Gamma_n$ on the circle satisfying similar dynamical properties.
(Actually, $\Gamma_n $ embeds into
$\widetilde{\mathrm{PSL}}(2,\mathbb{R})$ for $n \!\geq\! 2$,
and its image in $\mathrm{PSL}(2,\mathbb{R})$ is isomorphic
to the so-called {\em Hecke groups}
$\langle u,v \!: u^2 \!=\!v^{n+1}  \rangle$ \cite[Chapter II, Example 28]{harpe}.)
Although this idea is new in the context, it is very natural. Indeed, the property to be
shown implies that $\Gamma_n$ is left-orderable, and left-orderability for countable groups
is equivalent to the existence of faithful actions by (orientation-preserving)
homeomorphisms of the real-line \cite{ghys,order}. In the present case,
the action on the circle appears by taking the quotient with respect
to the central cyclic subgroup $\langle a^{n+1} \rangle$.


\section{$\Gamma_n$ as a group of fractions}
\label{ore}

\hspace{0.45cm} Unless otherwise stated, in 
what follows we will only consider the case $n > 1$:
the case $n=1$ is elementary and we leave it to the reader.

We begin by noticing that $\Delta = a^{n+1}$
belongs to the center of $\Gamma_n$.\footnote{Although it will be not used in
this work, it is worth mentioning that the center of $\Gamma_n$ coincides with
the cyclic group generated by $a^{n+1}\!.$ This is a direct consequence of
\cite{picantin}. More elementary, this can be easily deduced by looking
at the embedding of $\Gamma_n$ in $\widetilde{\mathrm{PSL}}(2,\mathbb{R})$
to be discussed in \S \ref{positive}.} Indeed,
$$b \Delta = b a^{n+1} = (b a^n) a = (a b^{-1}) a
= a (b^{-1} a) = a (a^n b) = a^{n+1} b = \Delta b,
\qquad a \Delta = a^{n+2} = \Delta a.$$
A word in (positive powers of) $a,b$ (resp. $a^{-1},b^{-1}$) will be said
to be {\em positive} (resp. {\em negative}). It is {\em non-positive} (resp.
{\em non-negative}) if it is either trivial or negative (either trivial or
positive).

\vsp

\begin{prop} {\em Every element $w \!\in\! \Gamma_n$ may be written in the form \esp
$\bar{u} \Delta^{\ell}$ \esp for some non-negative word $\bar{u}$ and $\ell \in \mathbb{Z}$.}
\label{astuce}
\end{prop}

\noindent{\bf Proof.} In any word representing $w$, we may 
rewrite the negative powers of $a$ and $b$ using the relations 
\esp $a^{-1} = a^n \Delta^{-1}$ \esp and $b^{-1} = \Delta^{-1} a^n b a^n.$ 
\esp Since $\Delta$ belongs to the center of $\Gamma_n$, 
this shows the proposition.\footnote{This argument
is motivated by the fact that the presentation
\esp $\Gamma_n = \big\langle a,c \!: \esp c a c = a^n \big\rangle$ \esp
endows $\Gamma_n$ with a structure of a
{\em Garside group}: see \cite[Page 268]{deh}
and \cite[Example 2]{paris}. Indeed, as is well-known, Garside
groups are groups of fractions of the corresponding monoids.}
$\hfill\square$

\vspace{0.5cm}

Let us take a more careful look at the positive words in $a, b$.
Using the relation $b a^n b = a$, one easily concludes that
every $v \!\in\! \langle a,b \rangle^+$ may be written in the form
$$v = b^{n_0} a^{m_1} b^{n_1} \cdots b^{n_{k-1}} a^{m_k} b^{n_k} \Delta^{\ell},$$
where $n_i > 0$ for $i \in \{1,\ldots,k\}$, $n_0 \geq 0$, $n_k \geq 0$,
$m_i \!\in\! \{1,\ldots,n\}$, and $\ell \geq 0$. Moreover, if $m_i \!=\! n$
and $0 \!<\! i \!<\! k$, then we may replace \esp
$b^{n_{i-1}} a^{m_i} b^{n_i} = b^{n_{i-1} - 1} (b a^n b) b^{n_i -1}$
\esp by \esp $b^{n_{i-1} - 1} a b^{n_i -1}$. \esp
This is also possible for $i\!=\!0$ (resp. $i\!=\!k$) when $n_0 \!>\! 0$
(resp. $n_k \!>\! 0$). Performing these reductions as far as possible,
we conclude that $v$ may be written in the form
$$v = b^{n_0} a^{m_1} b^{n_1} \cdots b^{n_{k-1}} a^{m_k} b^{n_k} \Delta^{\ell},$$
so that the following properties are satisfied:

\vsp

\noindent-- (i) \esp $n_i > 0$ for $0<i<k$, $n_0 \geq 0$, $n_k \geq 0$.

\vsp

\noindent-- (ii) \esp $m_i \!\in\! \{1,\ldots,n-1\}$ for $1 \!<\! i \!<\!k$.

\vsp

\noindent-- (iii) \esp $m_1$ lies in $\{1,\ldots,n-1\}$ (resp. 
$\{1,\ldots,n\}$) if $n_0 > 0$ (resp. $n_0 = 0$); \esp \esp similarly, $m_k$ 
lies in $\{1,\ldots,n-1\}$ (resp. $\{1,\ldots,n\}$) if $n_k > 0$ (resp. $n_k = 0$).

\vsp

\noindent-- (iv) \esp $\ell \geq 0$.

\vsp

\noindent Here, for $k=0$, an expression as above should be understood as 
$b^{n_0} \Delta^{\ell}$, where $n_0 \geq 0$.

\vsp\vsp

Therefore, by Proposition \ref{astuce}, every element $w \in \Gamma_n$
may be written in the form
\begin{equation}\label{p}
w \esp
= \esp b^{n_0} a^{m_1} b^{n_1} \cdots b^{n_{k-1}} a^{m_k} b^{n_k} \Delta^{\ell} \esp
= \esp u \Delta^{\ell},
\end{equation}
where properties (i), (ii), and (iii) above are satisfied, and $\ell \in \mathbb{Z}$.
Such an expression will be said to be a {\em normal form} for $w$.


\section{Eliminating numerators or denominators}
\label{en-uno-solo}

\hspace{0.45cm} Let $w = u \Delta^{\ell}$ be a normal form of a nontrivial
element $w \!\in\! \Gamma_n$. Our task consists in showing that $w$ is
either positive or negative in $a,b$. If $u$ is trivial, then $w$
is positive or negative according to the sign of $\ell$. If $u$ is
nontrivial and $\ell \geq 0$, then $w$ is positive. Assume throughout
that $u$ is nontrivial and $\ell < 0$. We will show that, in this
situation, $w$ is negative.

\vspace{0.35cm}

\noindent{\bf Case I.} We have $u = b^r$ for some positive integer $r$.

\vsp\vsp

In this case, the relation $b a^n b = a$ yields $a b a^{-1} = a^{-(n-1)} b^{-1}$,
hence
$$w \esp = \esp b^r a^{-1} a^{-n} \Delta^{\ell-1} \esp = \esp
a^{-1} (aba^{-1})^r a^{-n} \Delta^{\ell - 1} \esp = \esp
a^{-1} (a^{-(n-1)} b^{-1})^{r} a^{-n} \Delta^{\ell - 1}.$$

\vspace{0.35cm}

\noindent{\bf Case II.} The element $u$ does not belong
to $\langle b \rangle$:

\vsp\vsp

Let us consider the normal form (\ref{p}). There are two possibilities:

\vsp

\noindent (i) If $n_k = 0$, then using the relation
$a b a^{-1} = a^{-(n-1)} b^{-1}$ we obtain
\begin{eqnarray*}
w
&=& b^{n_0} a^{m_1} b^{n_1} \cdots a^{m_{k-1}} b^{n_{k-1}} a^{m_k} \Delta^{\ell}\\
&=& b^{n_0} a^{m_1} b^{n_1} \cdots a^{m_{k-1}-1}
    \underline{a b^{n_{k-1}} a}^{\!-1} a^{m_k - n} \Delta^{\ell-1}\\
&=& b^{n_0} a^{m_1} b^{n_1} \cdots a^{m_{k-1} - 1} (a^{-(n-1)} b^{-1})^{n_{k-1}}
    a^{m_k - n} \Delta^{\ell - 1}\\
&=& b^{n_0} a^{m_1} b^{n_1} \cdots a^{m_{k-2}} b^{n_{k-2}} a^{m_{k-1} - n}
    b^{-1} (a^{-(n-1)} b^{-1})^{n_{k-1} - 1} a^{m_k - n} \Delta^{\ell - 1}\\
&=& b^{n_0} a^{m_1} b^{n_1} \cdots a^{m_{k-2}- 1} \underline{a b^{n_{k-2}} a}^{-1} \esp
    a^{m_{k-1} - n + 1} b^{-1} (a^{-(n-1)} b^{-1})^{n_{k-1} - 1} a^{m_k - n} \Delta^{\ell - 1}.
\end{eqnarray*}

\vsp

\noindent (ii) If $n_k > 0$, then
\begin{eqnarray*}
w
&=& b^{n_0} a^{m_1} b^{n_1} \cdots a^{m_{k}-1} \underline{a b^{n_k} a}^{\!-1} a^{-n}
    \Delta^{\ell-1}\\
&=& b^{n_0} a^{m_1} b^{n_1} \cdots a^{m_{k} - 1} (a^{-(n-1)} b^{-1})^{n_{k}} a^{-n}
    \Delta^{\ell - 1}\\
&=& b^{n_0} a^{m_1} b^{n_1} \cdots a^{m_{k} - n}
    b^{-1} (a^{-(n-1)} b^{-1})^{n_{k} - 1} a^{-n} \Delta^{\ell - 1}\\
&=& b^{n_0} a^{m_1} b^{n_1} \cdots a^{m_{k-1}-1} \underline{ab^{n_{k-1}} a}^{\!-1}
    \esp a^{m_{k} - n+1} b^{-1} (a^{-(n-1)} b^{-1})^{n_{k} - 1} a^{-n} \Delta^{\ell - 1}.
\end{eqnarray*}

\vsp

\noindent The main point here is that in (i) we have
$m_k - n \leq 0$ and $m_{k-1} - n + 1 \leq 0$. Similarly,
in (ii) we have $m_k - n + 1 \leq 0$. This allows repeating
the argument. Proceeding in this way as far as possible, it
is easy to see that the final output will be a negative
word representing $w$.


\section{No positive word represents the identity}
\label{positive}

\hspace{0.45cm} We begin with a proof that applies
to $B_3 = \langle a,b\!: b a^2 b = a \rangle$. The
case of $\Gamma_n$ will be treated with a similar idea.

\vspace{0.15cm}

\begin{prop} {\em No element in $\langle a,b \rangle^+ \subset B_3$ represents the identity.}
\end{prop}

\noindent{\bf Proof.} Consider the representation
of $B_3$ in $\mathrm{PSL}(2,\mathbb{R})$ given by
$$a \to \bar{a} =
\left[
\begin{array}
{cc}
0 & 1  \\
-1 & 1  \\
\end{array}
\right],\qquad
b \to \bar{b} =
\left[
\begin{array}
{cc}
1 & 0  \\
1 & 1  \\
\end{array}
\right].$$
Denote by $U$ (resp. $V$) the projection of $\{(x,y) \!: x > y > 0\}$
(resp. $\{(x,y)\!: 0<x<y \}$) into $\mathbb{P}^1 (\mathbb{R})$.
A direct computation shows that
$$\bar{a}
\left[
\begin{array}
{cc}
x  \\
y \\
\end{array}
\right]
=
\left[
\begin{array}
{cc}
y   \\
y-x \\
\end{array}
\right],
\quad\quad
\bar{b}^n
\left[
\begin{array}
{cc}
x  \\
y \\
\end{array}
\right]
=
\left[
\begin{array}
{cc}
x      \\
nx + y \\
\end{array}
\right],$$
which easily yields $\bar{a}(V) \subset U$ and
$\bar{b}^n (U \cup V) \subset V$, for all $n > 0$. (See Figure 3 below.)

Now given an element $w \in \langle a,b \rangle^+$, let us write it in the form
$$w = b^{n_0} a^{m_1} b^{n_1} \cdots a^{m_k} b^{n_k} a^{3r},$$
where $k\geq 0$, $n_0,n_k,r$ are non-negative, $n_i > 0$ for the other indexes $i$,
and $m_i \!\in\! \{1,2\}$, with $m_i = 1$ for $1 < i < k$, and $m_1 = 1$
(resp. $m_k = 1$) when $n_0 > 0$ (resp. $n_k > 0$).
Notice that $\bar{a}^3 = id$. Assume that $w$ is not a power
of $a^3$. In this case, to show that $w \neq id$, it suffices to prove that
$$\bar{w} = \bar{b}^{n_0} \bar{a}^{m_1} \bar{b}^{n_1} \cdots \bar{a}^{m_k} \bar{b}^{n_k}$$
is nontrivial in $\mathrm{PSL}(2,\mathbb{R})$. Now using the relation
$\bar{b} \bar{a}^2 \bar{b} = \bar{a}$, one can easily check that,
unless $\bar{w}$ is conjugate to a power of $ab = \sigma_1$, it is conjugate to
a word $\bar{w}'$ in $\bar{a}$ (with no use of $\bar{a}^2$) and positive powers of
$\bar{b}$ which either begins and finishes with a power of $\bar{b}$, or begins and
finishes with $\bar{a}$. Since $\sigma_1$ is not a torsion element, $\bar{w} \neq id$
when $\bar{w}$ is conjugate to $\sigma_1$. Otherwise, a ping-pong type argument shows
that either $\bar{w}'(U) \subset V$ or $\bar{w}'(V) \subset U$, hence $\bar{w}' \neq id$.
$\hfill\square$

\vspace{0.5cm}


\beginpicture

\setcoordinatesystem units <1cm,1cm>

\circulararc 15 degrees from 2.15 0.7
center at -5 1.6

\plot 2.15 0.7 2.07 0.87 /
\plot 2.15 0.7 2.25 0.85 /

\circulararc -15.7 degrees from 1.28 0.7
center at 8 1.6

\plot 1.3 2.55 1.35 2.35 /
\plot 1.3 2.55 1.15 2.35 /

\plot 3 3 -3 -3 /
\plot -3 3 3 -3 /
\plot 0 -3.5 0 3.5 /
\plot -3.7 0 3.7 0 /

\put{$\bar{a}$} at 2 1.45
\put{$\bar{b}^n$} at 1 1.75

\put{} at -8.4 3.5
\put{Figure 3} at 0 -4

\begin{footnotesize}
\put{$\{ (x,y) \! : x > y > 0 \}$} at 2.1 0.3
\put{$\{ (x,y) \! : y > x > 0 \}$} at 1.4 3.2
\end{footnotesize}

\endpicture


\vspace{0.65cm}

The representation considered above is obtained via the well-known
identification of $B_3$ to $\widetilde{\mathrm{PSL}}(2,\mathbb{Z})$,
followed by the quotient by the center $\langle a^3 \rangle$. Indeed,
with respect to the system of generators $\{f = a, \esp h = b^{-1} a\}$,
the presentation of $B_3$ becomes \esp $\langle f,h \!: f^3 = h^2 \rangle$.

It turns out that $\Gamma_n$ also embeds into
$\widetilde{\mathrm{PSL}}(2,\mathbb{R})$. To see this, we first
rewrite the presentation of $\Gamma_n$ in terms of $f = a$ and $h = b^{-1}a$:
$$\Gamma_n \esp = \esp \big\langle f,h \!: f^{n+1} = h^2 \big\rangle.$$
This presentation shows that $\Gamma_n$ corresponds to a central extension of the Hecke group
$$H(n\!+\!1) = \big\langle \bar{f},\bar{h} \!: \bar{f}^{n+1} = \bar{h}^2 =id \big\rangle.$$
A concrete realization of $H(n\!+\!1)$ inside $\mathrm{PSL}(2,\mathbb{R})$ arises
when identifying $\bar{f}$ to the circle rotation of angle $\frac{2\pi}{n+1}$,
and $\bar{h}$ to the hyperbolic reflexion with respect to the geodesic joining
$p_n \!=\! \bar{f}^n(p)$ and $p=p_0$ for some point $p \!\in\! \clo$. This realization
allows embedding $\Gamma_n$ into $\widetilde{\mathrm{PSL}}(2,\mathbb{R})$ by identifying
$f \!\in\! \Gamma_n$ to the lifting of $\bar{f}$ to the real line given by
$x \mapsto x + \frac{2\pi}{n+1}$, and $h$ to the (unique) lifting $h$ of
$\bar{h}$ satisfying \esp $x \leq h(x) \leq x + 2\pi$
\esp for all $x \in \mathbb{R}$.\footnote{Actually, the arguments
given so far only show that the above identifications induce a group
homomorphism from $\Gamma_n$ into $\widetilde{\mathrm{PSL}}(2,\mathbb{R})$,
and the injectivity follows from the arguments given below combined with
the result of \S \ref{en-uno-solo}.}

\vspace{0.15cm}

The dynamics of the action of $H(n\!+\!1)$ on the circle is illustrated in Figure 4.
Here, \esp $\bar{g} \!=\! \bar{f} \bar{h}^{-1} \!=\! \bar{f} \bar{h}$ \esp
is a parabolic M\"obius transformation fixing $p_0$ and sending
$p_n$ into $p_1$, where \esp $p_i \!=\! \bar{f}^i (p)$ \esp for
$0 \!\leq\! i \!\leq\! n$. Using this action, we now proceed to show that no
element $w$ in $\langle a,b \rangle^+ \subset \Gamma_n$ represents the identity.

\vspace{0.5cm}


\beginpicture

\setcoordinatesystem units <0.9cm,0.9cm>

\circulararc 360 degrees from 3 0
center at 0 0

\circulararc -269 degrees from 1.8 1.75
center at 0 0

\circulararc 175 degrees from 1 3
center at 0 3
\plot -1 3 -0.86 3.3 /
\plot -1 3 -1.12 3.3 /
\put{$\bar{f}$} at 0 4.3

\circulararc 175 degrees from 2.6 1.8
center at 1.9 2.34
\plot 1.25 2.9 1.5 3 /
\plot 1.25 2.9 1.4 3.2 /
\put{$\bar{f}$} at 2.7 3.3

\circulararc -175 degrees from -2.6 1.8
center at -1.9 2.34
\plot -2.6 1.8 -2.88 1.95 /
\plot -2.6 1.8 -2.58 2.05 /
\put{$\bar{f}$} at -2.6 3.3

\plot -1.8 1.75 -1.9 1.37 /
\plot -1.8 1.75 -2.18 1.55 /

\put{$\bar{g}$} at 0 -2.25

\put{$p = p_0$} at 0.1 3.2
\put{$p_1$} at -2.05 2.55
\put{$p_2$} at -3.1 1.1
\put{$p_n$} at 2.1 2.6
\put{$p_{n-1}$} at 3.4 1.1
\put{$p_{n-2}$} at 3.4 -1.1
\put{$p_3$} at -3.1 -1.1
\put{$p_{n-3}$} at 2.4 -2.6

\put{$\bullet$} at 0 3
\put{$\bullet$} at -2 2.2
\put{$\bullet$} at -2.9 0.75
\put{$\bullet$} at 2 2.2
\put{$\bullet$} at 2.9 0.75
\put{$\bullet$} at -2.9 -0.75
\put{$\bullet$} at 2.9 -0.75
\put{$\bullet$} at 2 -2.24

\put{} at -9.2 3.5
\put{Figure 4} at 0 -3.5

\endpicture


\vspace{0.6cm}

We begin by writing $w$ in the form
$$w = b^{n_0} a^{m_1} b^{n_1} \cdots a^{m_k} b^{n_k} a^{(n+1)r},$$
where $k \geq 0$, $n_0,n_k,r$ are non-negative, $n_i > 0$ for $0 \!<\! i \!<\! k$,
and $m_i \!\in\! \{1,2,\ldots,n\}$, with $m_i \neq n$ for $1 < i < k$, and
$m_1 \neq n$ (resp. $m_k \neq n$) when $n_0 > 0$ (resp. $n_k > 0$).
If $w$ were equal to the identity, then
$$\bar{w} = \g^{n_0} \bar{f}^{m_1} \g^{n_1} \cdots \bar{f}^{m_k} \g^{n_k}$$
would act trivially on the circle. Assume that
$w$ is not a power of $a$. Then one easily checks that,
unless $\bar{w}$ is a power of $\f \g$, it is conjugate either to some
$\bar{w}' \in \langle \f, \g \rangle^+$ beginning and ending by $\f$
and so that all exponents of $a$ lie in $\{1,\ldots,n-1\}$,
or to some $\bar{w}'' \in \langle \f, \g \rangle^+$ begining and
ending with $g$ with all exponents of $a$
in $\{1,\ldots,n-1\}$. Now, an easy ping-pong type argument
shows that $\bar{w}' \big( ]p_0,p_1[ \big) \!\subset ]p_1,p_n[$ and
$\bar{w}'' \big( ]p_1,p_n[ \big) \!\subset ]p_0,p_1[$, and hence
$\bar{w}' \neq id$ and $\bar{w}'' \neq id$.
Thus, to conclude the proof, we need to check that
neither $a$ nor $\f \g$ are torsion elements.

That $\f \g$ is not torsion follows from that it sends $[p,p_n]$ into the
subinterval $[p_1,p_2]$, and hence no iterate of it can equal the identity.
Finally, to see that $a$ is not torsion, just notice that $a$ maps
to the translation by $\frac{2\pi}{n+1}$ in
$\widetilde{\mathrm{PSL}}(2,\mathbb{R})$, and
hence has infinite order.\footnote{This also follows
from the main result of \cite{torsion-deh}.}


\section{Dehornoy-like orderings}
\label{dehornoy}

\hspace{0.45cm} In what follows, we will denote by $\preceq_{n}$ the
left-ordering on $\Gamma_n$ whose positive cone is $\langle a,b \rangle^+$.
Using $\preceq_n$, we will define an analog of the Dehornoy ordering.

We begin by recalling the Dehornoy ordering on
$B_3$. Consider the Artin (standard) presentation
$$B_3 = \big\langle \sigma_1, \sigma_2 \!: \esp
\sigma_1 \sigma_{2} \sigma_1 = \sigma_{2} \sigma_1 \sigma_{2} \big\rangle.$$
Following Dehornoy \cite{dehornoy-libro}, an element of $B_3$ is said to be
$1$-positive if it may be written as a word of the form
$$\sigma_2^{n_0} \sigma_1 \sigma_2^{n_1} \sigma_1 \cdots
\sigma_2^{n_{k-1}} \sigma_1 \sigma_2^{n_k},$$
where $n_i \!\in\! \mathbb{Z}$.
It is said $2$-positive if it is of the form $\sigma_2^n$ for some $n \!>\! 0$. An element in
$B_3$ is said to be $D$-positive if it is either $1$-positive or $2$-positive. The remarkable
result of Dehornoy (in the case of $B_3$) asserts that the set of $D$-positive elements is the
positive cone of a left-ordering $\preceq_D$ on $B_3$. The proof given by Dehornoy as well as
many subsequent proofs are very intricate (see \cite{DDRW} for a detailed discussion on this).
Nevertheless, a short proof using the ordering $\preceq_2$ may be given. What follows is
inspired from \cite{order} (see Examples 3.35 and 3.36 therein).

Before continuing our discussion, recall that a subgroup $\Gamma_0$ of a
left-ordered group $(\Gamma,\preceq)$ is said to be {\em $\preceq$-convex} if $g$ 
belongs to $\Gamma_0$ whenever $h_1 \prec g \prec h_2$ for some $h_1,h_2$ in $\Gamma_0$.
Convex subgroups are very useful for defining new orders: If $(\Gamma,\preceq)$ and
$\Gamma_0$ are as above and $\preceq'$ is any left-ordering on $\Gamma_0$, then the
{\em extension} of $\preceq'$ by $\preceq$ is the left-ordering on $\Gamma$ whose
positive cone is
$$P^+ = P^+_{\preceq'} \cup \big( P^+_{\preceq} \setminus \Gamma_0 \big).$$
Finally, we denote by $\bar{\preceq}$ the {\em reverse ordering} of $\preceq$,
that is, the left-ordering defined by $g \esp \bar{\succ} \esp id$ if
and only if $g \prec id$.

\vspace{0.15cm}

\begin{lem} \label{unicos-convexos}
{\em For each $n \!\in\! \mathbb{N}$, the subgroup $\langle b \rangle
\! \subset \! \Gamma_n$ is $\preceq_n$-convex. Moreover, the only $\preceq_n$-convex
subgroups of $\Gamma_n$ are $\{ id \}$, $\langle b \rangle$, and $\Gamma_n$ itself.}
\end{lem}

\noindent{\bf Proof.} Let $c \in \Gamma_n$ be such that $b^r \prec_n c \prec_n b^s$. Assume
that $c$ is $\preceq_n$-positive (the other case is analogous). If $c$ does not belong
to $\langle b \rangle$, then it may be written in the form $w_1 a w_2$, where $w_1$ and
$w_2$ are (perhaps empty) words on non-negative powers of $a$ and $b$. The inequality
$c \prec_n b^s$ yields $w = b^{-s} w_1 a w_2 \prec_n id$. Introducing the identity
$a = ba^2 b$ several times, one easily shows that $w$ may be rewritten as
$w = w_1' w_2$, where $w_1'$ only uses positive powers of $a$ and $b$.
Thus, $w$ is $\preceq_n$-positive, which is a contradiction.

To show that the only $\preceq_{n}$-convex subgroups of $\Gamma_n$ are $\{ id \}$,
$\langle b \rangle$, and $\Gamma_n$ itself, we proceed by contradiction. Clearly,
$\langle b \rangle$ does not contain any nontrivial convex subgroup. Suppose
that there exists a $\preceq_{n}$-convex subgroup $N$ of $\Gamma_n$ such that
\esp $\langle b \rangle \subsetneq N \subsetneq \Gamma_n$. \esp Let $\preceq'$,
$\preceq''$, and $\preceq'''$, be the left-orderings defined on
$\langle b \rangle$, $N$, and $\Gamma_n$, respectively, by:

\vspace{0.1cm}

\noindent -- $\preceq'$ is the restriction of $\preceq_{n}$ to $\langle b \rangle$,

\vspace{0.1cm}

\noindent -- $\preceq''$ is the extension of $\preceq'$ by the restriction of
$\bar{\preceq}_{n}$ to $N$,

\vspace{0.1cm}

\noindent -- $\preceq'''$ is the extension of $\preceq''$ by $\preceq_{n}$.

\vspace{0.1cm}

\noindent The order $\preceq'''$ is different from $\preceq_{n}$ (the $\preceq_{n}$-negative
elements in $N \setminus \langle b \rangle$ are $\preceq'''$-positive), but its positive cone
still contains the elements \esp $a,b$. \esp Nevertheless, this is impossible, since these
elements generate the positive cone of $\preceq_{n}$. $\hfill\square$

\vspace{0.5cm}

Now let $\bar{\preceq}_n$ be the reverse ordering of $\preceq_n$, and let $\preceq_n'$ be 
the ordering of $\Gamma_n$ obtained as the extension of $\bar{\preceq}_n$ (restricted to
$\langle b \rangle$) by $\preceq_n$. We claim that, for $n \!=\! 2$ ({\em i.e.} for 
$B_3$), $\preceq_n'$ coincides with the Dehornoy ordering $\preceq_D$. Indeed, if 
$c \in \langle b \rangle$ is $\preceq_2'$-positive, then it is a negative 
power of $b = \sigma_2^{-1}$, hence $\preceq_D$-positive. If 
$c \in B_3 \setminus \langle b \rangle$ is $\preceq'_2$-positive, then 
it may be written as a word using only positive powers of $a$ and $b$. Replacing 
$a = \sigma_1 \sigma_2$ and $b=\sigma_2^{-1}$, this allows writing $c$ as a word 
where only positive powers of $\sigma_1$ are used. In particular, $c$ is 
$\preceq_D$-positive. We thus conclude that the positive cone of $\preceq_{D}$ 
contains that of $\preceq_2'$. Conversely, if a nontrivial element $c$ is not 
$\preceq_2'$-positive, then $c^{-1}$ is $\preceq_2'$-positive, hence 
$\preceq_D$-positive; thus, $c$ is neither $\preceq_D$-positive. This 
shows that the positive cone of $\preceq_2'$ contains that of $\preceq_D$.


\vspace{0.35cm}

The equivalence between $\preceq_D$ and $\preceq_2'$ gives a new proof of Dehornoy's
theorem (for $B_3$). It also motivates the following definition.

\vsp

\begin{defn} For each $n \geq 2$ the left-ordering $\preceq_n'$ on
$\Gamma_n$ will be called the {\em Dehornoy-like ordering} of $\Gamma_n$.
\end{defn}

\vsp

As in the case of $B_3$, an element $c \!\in\! \Gamma_n$ is $\preceq'_n$-positive if
either $c = b^{-k}$ for some $k \geq 1$, or it may be written in the form
$$c = b^{n_0} a b^{n_1} a \cdots b^{n_{k-1}} a b^{n_k}$$
for some $n_i \! \in \! \mathbb{Z}$ (with $k \geq 1$). Notice that the smallest
positive element of $\preceq'_n$ is $b^{-1}$. Moreover, the family of
$\preceq_n'$-convex subgroups of $\Gamma_n$ coincides with that of
$\preceq_n$-convex ones, that is,
$\{ id \}, \langle b \rangle, \Gamma_n$ (see \cite[Remark 3.34]{order}).
The following proposition (and its proof) extends \cite[Theorem D]{order}.

\vspace{0.15cm}

\begin{prop} {\em The positive cone of the Dehornoy-like ordering
$\preceq_n'$ of $\Gamma_n$ is not finitely generated as a semigroup.}
\end{prop}

\noindent{\bf Proof.} Following \cite[Example 8.2]{braids}, we will show that the
sequence of conjugates $b^k a (\preceq'_n)$ converges to $\preceq'_n$ in a nontrivial way.
Here, $b^k a (\preceq'_n)$ is the left-ordering whose positive cone is the conjugate \esp
$b^k a P_{\preceq'_n} (b^k a)^{-1}$
of $P_{\preceq_n'}$. \esp Saying that $b^k a(\preceq'_n)$ converges to
$\preceq'_n$ in a nontrivial way means that, though $b^k a (\preceq'_n)$ does not coincide
with $\preceq'_n$ for $k$ large enough, given finitely many $\preceq'_n$-positive elements
$c_1,\ldots,c_r$, these elements are also positive with respect to $b^k a (\preceq'_n)$ for
$k$ large enough. Such a convergence implies that the positive cone of $\preceq'_n$ cannot
be finitely generated. Indeed, if it were generated by $c_1,\ldots,c_r$, then these elements
would be positive for $b^ka(\preceq'_n)$ for $k$ large enough. This would imply that
$b^ka(\preceq'_n)$ coincides with $\preceq'_n$ for large $k$, which is a contradiction.

If $c_i$ does not belong to $\langle b \rangle$, then $c_i$ may be written in the form \esp
$c_i = b^{n_0} a \bar{w}$ \esp for some $n_0 \!\in\! \mathbb{Z}$ and a certain
$\bar{w}$ containing no negative power of $a$. We then have
$$(b^ka)^{-1} c_i b^k a = a^{-1} b^{-k + n_0} a \bar{w} b^k a.$$
For $k > n_0$, the relation $a^{-1}b^{-1}a = a^{n-1} b$ yields
$$(b^ka)^{-1} c_i b^k a = (a^{n-1} b)^{k - n_0} \bar{w} b^k a.$$
The right-side expression above contains only positive powers of $a$, thus showing that
$c_i$ is positive with respect to $b^k a (\preceq'_n)$ provided that $k > n_0$.

If $c_i$ belongs to $\langle b \rangle$, then $c_i = b^{-r}$ for some $r \in \mathbb{N}$.
This yields
$$(b^ka)^{-1} c_i b^k a = a^{-1} b^{-k} b^{-r} b^k a = a^{-1} b^{-r} a = (a^{n-1}b)^r.$$
The right-side expression contains only positive powers of $a$,
hence it is $\preceq'_n$-positive.

Finally, to show that $b^ka(\preceq'_n)$ and $\preceq'_n$ do not coincide, it suffices
to notice that the smallest positive element of the former ordering, namely
$(b^ka)^{-1} b^{-1} b^k a = a^{-1} b^{-1} a$, is different from $b^{-1}$,
which is the smallest positive element of $\preceq'_n$. $\hfill\square$

\vspace{0.12cm}

\begin{rem} The very same argument of the proof above shows that
$b^k a (\preceq_n)$ also converges to $\preceq_n'$ as $k$
goes to infinity.
\end{rem}

\vspace{0.07cm}

Another relevant property of the Dehornoy ordering on $B_3$ is the so-called
{\em Property S}: All conjugates of $\sigma_1$ and $\sigma_2$ are $\preceq_D$-positive.
We were not able to reprove this property with our methods. More importantly, we do not
know whether an analog of this property holds for all Dehornoy-like orderings.


\section{Some questions and comments}
\label{questions}

\hspace{0.45cm} One may address plenty of questions on the structure of the groups
$\Gamma_n$. However, we would like to focus on certain aspects related to group orderability.

\vspace{0.5cm}

\noindent{\bf $C$-orderability and local indicability.}
Let us recall that a group is said to be {\em $C$-orderable} if it admits a left-ordering
$\preceq$ satisfying $fg^k \! \succ \! g$ for all $f,g$ positive and all $k \geq 2$ (see
\cite{order}). Such an ordering is said to be {\em Conradian}. A remarkable theorem
of Brodskii \cite{brodski}
asserts that torsion-free, 1-relator groups are $C$-orderable. Indeed, such a group is
necessarily {\em locally indicable} \cite{brodski,howie} (that is,
each of its finitely generated subgroups
surjects into $\mathbb{Z}$), and local indicability is equivalent to $C$-orderability
(see \cite[\S 3]{order} for a discussion on this point).

Now notice that, since the groups $\Gamma_n$ are left-orderable, they are torsion-free.
(This also follows from \cite{torsion-deh}.) By the discussion above, they are
$C$-orderable.\footnote{Notice that the $C$-orderability of $\Gamma_n$ does not
follows from that it embeds into $\widetilde{\mathrm{PSL}}(2,\mathbb{R})$.
Indeed, $\widetilde{\mathrm{PSL}}(2,\mathbb{R})$ contains finitely generated 
groups with trivial first cohomology, as for example the lifting of the 
$(2,3,7)$-triangle group \cite{bergman,Th}.} For example, the local indicability 
of $B_3 \!\sim\! \Gamma_2$ comes from the well-known exact sequence
$$0 \esp \esp \longrightarrow \esp\esp [B_3,B_3] \sim F_2 \esp\esp \longrightarrow
\esp\esp B_3 \esp\esp \longrightarrow \esp\esp B_3 / [B_3,B_3] \sim \mathbb{Z}
\esp\esp \longrightarrow \esp\esp 0$$
and the fact that free groups are locally indicable (this last result goes back
to Magnus \cite{MKS}). We point out, however, that the orderings $\preceq_n$ and
$\preceq_n'$ are not Conradian (for $n \!>\! 1$):

\vsp

\noindent -- For $\preceq_n$, notice that $a \!\succ_n\! id$ and
$b \!\succ_n\! id$, though \esp $a^{-1}ba^n \!=\! b^{-1} \!\prec_n\! id$, \esp thus
$b a^n \!\prec_n\! a$.

\vsp

\noindent -- For $\preceq_n'$, we have $ab^{2} \succ_n' \! id$
and $ab \succ_n' \! id$. Now from \esp
$a^{-1} b a = b^{-1} a^{-(n-1)}$ \esp we obtain
\begin{eqnarray*}
(ab)^{-2} (ab^{2}) (ab)^{4}
&=& b^{-1} \underline{a^{-1} b a} b (ab)^3 =
    b^{-1} b^{-1} a^{-(n-1)} b (ab)^3\\
&=& b^{-2} a^{-(n-2)} \underline{a^{-1} b a} b (ab)^2\\
&\vdots&\\
&=& b^{-2} a^{-(n-2)} b^{-1} a^{-(n-2)} b^{-1}
a^{-(n-2)} b^{-1} a^{-(n-1)} b \esp\esp\esp\esp
\prec_n' \esp\esp\esp\esp id,
\end{eqnarray*}
hence \esp $a b^{2}  \big( (ab)^2 \big)^{2} \prec_n' (ab)^2.$

\vsp

\noindent As a more sophisticated argument, let us mention Conradian orderings 
with finitely many convex groups ({\em c.f.} Lemma \ref{unicos-convexos}) may 
only exist on solvable groups (see for instance \cite[\S 1.3]{F}), and the
groups $\Gamma_n$ (with $n \!>\! 1$) are non-amenable (to see this, just 
notice that the actions on the circle constructed in \S \ref{positive} 
have no invariant probability measure).

\vspace{0.5cm}

\noindent{\bf Other positive cones generated by two elements.} \hspace{0.015cm} 
It is interesting to compare the groups $\Gamma_n$ with the Baumslag-Solitar groups
\esp $BS_{1,n} \esp = \esp \big\langle a,b\!: \esp b^{-1} a^n b =a \big\rangle.$ \esp
Indeed, $BS_{1,n}$ is locally indicable (hence $C$-orderable), admits uncountably
many left-orderings, but only four $C$-orderings (all of which are bi-invariant).
Actually, this is nearly a characterization of these groups (see \cite{rivas}).
This gives some ``evidence'' for a positive answer to the following

\vspace{0.3cm}

\noindent{\bf Main Question.} Let $\Gamma$ be a group admitting a left-ordering 
whose positive cone is generated by (no more than) two elements. Is $\Gamma$ 
isomorphic to either $\mathbb{Z}$ or $\Gamma_n$ for some $n \geq 1$~?

\vspace{0.3cm}

Notice that the group \esp $\Gamma_{m,n} = \langle a,b\!: b a^n b = a^m\rangle$ \esp
is isomorphic to $\Gamma_{m+n-1}$ for all positive $m,n$, thus it belongs to the
family above. Indeed, the relator of $\Gamma_{m,n}$ may be written as
\esp $(b a^{n-1})^{-1} a^{m+n-1} (b a^{n-1})^{-1} = a$.

\vspace{0.5cm}

\noindent{\bf Positive cones generated by $k \! > \! 2$ elements.} According to
\cite{dub}, for each $n\geq 1$, the braid group $B_n$ admits a left-ordering
whose positive cone is generated by $n\!-\!1$ elements, namely
$$\sigma_{1} \sigma_{2} \cdots \sigma_{n-1}, \esp
(\sigma_{2} \sigma_3 \cdots \sigma_{n-1})^{-1}, \esp
\sigma_3 \sigma_4 \cdots \sigma_{n-1}, \esp \ldots \esp
, \esp (\sigma_{n-1})^{(-1)^n}.$$
Once again, the proof of this fact given in \cite{dub} uses Dehornoy's theory. We
were not able to extend our approach to simplify and/or generalize this phenomenon.
One of the difficulties lies in that, with the generators above, the natural
presentations of $B_n$ are not Garside. We expect, however, that some alternative
approach should yield an answer for the following

\vspace{0.3cm}

\noindent{\bf Main Problem.} For each $k > 3$, find an infinite family of groups
(including both $B_{k-1}$ and the Tararin groups $T_k$ from \cite[\S 4.2]{rivas})
all of which admit left-orderings with a positive cone generated (as a 
semigroup) by $k$ elements.


\vspace{0.8cm}

\noindent{\bf \Large{Acknowledgments.}} It is a pleasure to thank P. Dehornoy for 
useful references on Garside groups as well as many encouragements, \'E. Ghys for 
a clever suggestion, A. Glass for comments and corrections, and B. Wiest for 
explanations on the geometry of braid groups.

This work was funded by the PBCT/Conicyt Research Network on Low Dimensional Dynamics 
and Fondecyt Project 1100536.


\vspace{0.2cm}

\begin{small}


\vspace{0.3cm}

\noindent Andr\'es Navas\\

\noindent Dep. de Matem\'aticas,
Fac. de Ciencia, Univ. de Santiago\\

\noindent Alameda 3363, Estaci\'on Central, Santiago, Chile\\

\noindent E-mail address: andres.navas@usach.cl\\

\end{small}

\end{document}